\documentclass{amsart}
\usepackage{amsthm}
\usepackage{mathrsfs}

\newtheorem{defi}{Definition}
\newtheorem{thm}{Theorem}
\newtheorem{conj}{Conjecture}

\def\N{\mathbb{N}}
\def\Z{\mathbb{Z}}
\def\Q{\mathbb{Q}}
\def\U{\mathcal{U}}
\def\ng{\mathbb{N}[G]}
\def\nZ{\mathbb{N}[\Z/2\Z]}
\def\comp{\mathrm{Comp}\left(\ng\right)}
\def\pr{\mathrm{Pr}\left(\ng\right)}

\begin{document}
\title{Composites In Semirings of Boolean Groups}

\begin{abstract}
We estimate the number of composite elements in the $n$-th grade of the group semiring of finite boolean groups. In view of this result
we may conjecture that the composites in the semiring of finite groups are thinly dispersed.
\end{abstract}

\author{Kamalakshya Mahatab} 
\address{Kamalakshya Mahatab, Department of Mathematical Sciences, NTNU Trondheim, Norway}
\email{accessing.infinity@gmail.com, kamalakshya.mahatab@ntnu.no}


\maketitle
\section{Introduction}
A semiring $S$ is a commutative monoid with an addition \lq$+$\rq \ and also a monoid with a multiplication \lq$.$\rq \ such that 
the multiplication is distributive over addition \cite{golan}. In this article, we will consider group semirings $S$ with an additive 
identity \lq$0$\rq \ and with a multiplicative identity \lq$1$\rq. We will define units in a semiring $S$ as follows:
\begin{defi}
 A non-zero element $u\in S$ is called a unit in $S$, if there exist an element $v\in S$ such that $u.v=1$. 
 We will denote the group of units of $S$ as $\U(s)$.  
\end{defi}
In particular $1\in \U(s)$. Now we may define reducible and irreducible elements in a group semiring.

\begin{defi}
We will define a non-zero non-unit element $w\in S$ as a prime, if for every $u, v\in S$ with $uv=w$, we have either $u$ or $v$ is a unit.

A non-zero non-unit element of $S$ which is not prime will be called a composite.
\end{defi}

We will consider the semirings constructed from finite groups.
\begin{defi}
Let
\[\mathbb N:=\{0, 1, 2, \ldots\},\]
and $G$ be a finite group. We define the group semiring $\ng$ as
\[\ng=\left\{\sum_{g\in G}a_g g: a_g\in \mathbb N\right\}\]
with point-wise addition and multiplication inherited from the group $G$:
\begin{align*}
\left(\sum_{g\in G} a_g g\right) + \left(\sum_{g\in G}b_gg\right) &= \sum_{g\in G}(a_g+b_g)g,\\
\left(\sum_{g\in G} a_g g\right) . \left(\sum_{g\in G}b_gg\right) &=\sum_{g\in G}\left(\sum_{\substack{g_1, g_2\\ g_1g_2=g}} a_{g_1}b_{g_2}\right)g,
\end{align*}
where $a_g, b_g\in \N$ for all $g\in G$.
\end{defi}
\noindent
We will write $1.g$ as $g$, and $\sum_{g\in A}0.g$ as $0$ for any subset $A$ of $G$. 
Also observe that $\ng$ has a graded structure. For a $k\in\N$, define the $k$th grade of $\ng$ as
\[\ng_k:=\left\{\sum_{g\in G}a_g g\in \ng \colon \sum_{g\in G}a_g=k\right\}.\]
With the above notations, the set of units of $\ng$ is
\[\U\left(\ng\right)=\{g:g\in G\}=\ng_1.\]
Moreover, the group semiring $\ng$ is graded:
\begin{align*}
\ng_k+\ng_l &\subseteq \ng_{k+l},\\ 
\ng_k\ng_l &\subseteq \ng_{kl}, \text{ for all } k, l\in\mathbb N.
\end{align*}
We will write $\comp_k$ for the set of composite elements of $\ng_k$ and $\pr_k$ for the set of primes of $\ng_k$. In this article we will attempt
to give a \lq good\rq \ upper bound for $|\comp_k|$. As an initial step, we make the following observations.
\begin{thm}\label{thm:basic_conditions}
Let $p$ be any prime number and let $e$ denote the identity element in a finite group $G$. Then the following
statements hold:
 \begin{enumerate}
  \item[(i)] All the elements in the $p$-th grade of $\ng$ are prime:
  \[\pr_p=\ng_p.\]
  In particular, 
  \[\left|\pr_p\right|=\dbinom{p+|G|-1}{p-1}.\]
  \item[(ii)] For $n\geq 2,$
  \[\comp_n=\bigcup_{\substack{d|n, \\ 1<d<n}}\ng_d\ng_{n/d}.\]
  \item[(iii)] For $n\geq 2$ and for any $h\in G$ with $h\neq e,$ the element
  \[e+(n-1)h\]
 of $\ng$ is a prime.
 \end{enumerate}
\end{thm}

The last statement in the above theorem says that each grade of $\ng$ has a prime in it. Also we know that $\ng_p$ has no composite elements when
$p$ is prime. We may ask: How many composite elements are there in $\ng_n$ when $n$ is not a prime number? This is a 
difficult question to answer in general.
In this article, we give some answer to this question for the case $G=(\Z/2\Z)^l$. Even then, we will see that the upper bounds of $\comp_k$ requires 
some non-trivial results from number theory.

We may observe that $|\N[(\Z/2\Z)^l]_n|=\binom{n+2^l-1}{2^l-1}\asymp_l n^{2^l-1}$. Let 
$$\Theta_l(n):=\left|\text{Comp}\left(\N[(\Z/2\Z)^l]\right)_n\right|.$$ 
By Theorem~\ref{thm:basic_conditions}, 
\begin{align*}
\Theta_l(n)&\leq|\N[(\Z/2\Z)^l]_n-|\{g_1+(n-1)g_2: g_1, g_2 \in \N[(\Z/2\Z)^l\}|\\
&\leq \binom{n+2^l-1}{2^l-1}-2^{2l}+2^{l}.
\end{align*}
But this is a very crude upper bound. To compute an upper bound for $\Theta_l(n)$ when $n$ does not have too many prime factors, we need the following result of Kevin Ford.
\begin{thm}[Ford\cite{ford1, ford2}]\label{thm:ford}
 For $m\leq n$ and $a_1, a_2, b_1, b_2 \in \N$, we have
 \[\left|\{a_1+jb_1: 1\leq j\leq m\}\{a_2+jb_2: 1\leq j\leq n\}\right|\asymp\frac{mn}{(\log m)^{\delta}(\log\log m)^{3/2}},\]
 where \footnote{We will fix the notation $\delta=1-\frac{1+\log\log 2}{\log2}$ throughout this article. }
 \[\delta=1-\frac{1+\log\log 2}{\log2}.\] 
\end{thm}


We will only give a sketch of the proof of Theorem~\ref{thm:ford} as it is an easy consequence of Ford's proof. In \cite{ford1, ford2}, Ford proved that
\[|\{j_1j_2: 1\leq j_1, j_2\leq n \}|\asymp\frac{n^2}{(\log n)^{\delta}(\log\log n)^{3/2}},\quad n\rightarrow \infty.\] 
This proof also implies
\begin{equation}\label{eq:arthset}
 |\{j_1j_2: 1\leq j_1\leq m, j_2\leq n \}|\asymp\frac{mn}{(\log n)^{\delta}(\log\log n)^{3/2}},\quad \text{ when } m\leq n,\ m\rightarrow \infty.
\end{equation}
Further note that 
\begin{align}\label{eq:arthset2}
\notag
 &\left|\{(a_1+j_1b_1)(a_2+j_2b_2): 1\leq j_1\leq m, 1\leq j_2\leq n\}\right|\\
 &=|\{b_1j_1b_2j_2: 1\leq j_1\leq m, j_2\leq n\}|+ O(m+n).
\end{align}
So Theorem~\ref{thm:ford} follows from (\ref{eq:arthset}) and (\ref{eq:arthset2}).

Using the above estimate, we prove the following theorem.
\begin{thm}\label{thm:comp_estimate}
 Let $k$ be a fixed positive integer. Let $P^-(n)$ denote the smallest prime factor of $n$ for a positive integer $n\geq 2$.
 Further, if $n$ has at most $k$ prime factors, then 
 \[\Theta_l(n)\ll_{k, l}\left(\frac{n}{(\log P^-(n))^{\delta}(\log\log P^-(n))^{3/2}}\right)^{2^l-1} \quad\quad \text{as} \quad P^-(n)\rightarrow\infty.\]
 For $l=1$,
 \[\Theta_1(n)\asymp_k\frac{n}{(\log P^-(n))^{\delta}(\log\log P^-(n))^{3/2}} \quad\quad \text{as} \quad P^-(n)\rightarrow\infty.\]
 \end{thm}
 
Theorem~\ref{thm:comp_estimate} suggests that for any finite group $G$, the number of composites is much smaller 
than the number of primes in the $n$-th grade of the group semiring. In other words,
\begin{conj}\label{conj}
Let $G$ be a finite group and let the number of prime factors of $n$ be bounded by $k$.
Then 
\[\lim_{P^-(n)\rightarrow\infty}\frac{\left|\comp_n\right|}{\left|\ng_n\right|}=0.\]
\end{conj}

We would like to mention that we could not find any relevant literature on the above conjecture. Similar questions can be formulated for primes in 
different classes of positive matrices (see \cite{PHS}). 

To prove Theorem~\ref{thm:comp_estimate}, we used the fact that the representations of $(\Z/2\Z)^l$ are one-dimensional and rational. In other words, $\Q\left[(\Z/2\Z)^l\right]$ splits
completely in one-dimensional representations, which allows us to rephrase our question on primes in group semirings in terms of certain distribution of prime numbers. However, this is not true
for other groups, which is the main obstacle in proving Conjecture~\ref{conj}.

\section{Proofs Of The Theorems}
\subsection{Proof of Theorem \ref{thm:basic_conditions}}
 Proofs of (i) and (ii) follows from the definitions. To prove (iii), let 
 $$\sum_{g\in G}a_g g, \ \sum_{g\in G}b_g g \ \in \ \nZ$$ 
 be such that
 \begin{align*}
  &\left(\sum_{g\in G}a_g g\right) \left( \sum_{g\in G}b_g g \right)= e+(n-1)h.
 \end{align*}
We have the following identities: 
\begin{align}\label{eq:ag-bg-1}
 &\sum_{g\in G}a_g b_{g^{-1}}=1,\\
 \label{eq:ag-bg-h}
 &\sum_{g\in G}a_g b_{g^{-1}h}=n-1,\\
 \label{eq:ag-bg-etc}
 &\sum_{g\in G}a_g b_{g^{-1}h'}=0 \ \text{ for } h'\neq e, h.
\end{align}
Let for some $g_0\in G, a_{g_0}\neq 0$ and it contributes to the sum in (\ref{eq:ag-bg-1}). Then 
\begin{equation}\label{eq:ag0-bg0}
a_{g_0}=b_{g_0^{-1}}=1. 
\end{equation}
From (\ref{eq:ag-bg-etc}), we get 
\begin{equation}\label{eq:bg0-inv-h}
 b_{g_0^{-1}h'} = 0 \ \text{ for } \ h'\neq e, h.
\end{equation}
We prove (iii), if we have $b_{g_0^{-1}h}=0$; so assume that it is non-zero.
Then by (\ref{eq:ag-bg-1})
\begin{align}\label{eq:ag0-h-inv}
 a_{g_0h^{-1}}=0.
\end{align}
From (\ref{eq:ag-bg-h}), (\ref{eq:ag0-bg0}) and (\ref{eq:bg0-inv-h}), we get
\begin{equation}
 b_{g_0^{-1}h} = n-1.
\end{equation}
Finally, we conclude that 
\[\sum_{g\in G}a_g g = g_0 \ \text{ and } \ \sum_{g\in G}b_g g = g_0^{-1} + (n-1)g_0^{-1}h. \]
This proves (iii).
\subsection{Proof of Theorem \ref{thm:comp_estimate}}
We will first prove the theorem for the case $l=1$. 
%
Let 
\[\Z/2\Z=\{\alpha_0, \alpha_1\},\]
with $\alpha_0$ being the identity element. Define a map $\Psi$ as follows:
\begin{align*}
&\Psi: \nZ\longrightarrow\Z &\\
&\quad a\alpha_0+b\alpha_1\longrightarrow a-b.&
\end{align*}
The map $\Psi$ is a semiring homomorphism; in other words, it preserves multiplication and addition of $\ng$ in the ring 
$\Z$. Let $\Psi_n$ be the restriction of $\Psi$ to $\nZ_n$. We may also observe that $\Psi_n$ is a one to one map.
If $n=m_1m_2$ for $m_1, m_2>1$, then
\begin{align*}
&\Psi_n\left(\nZ_{m_1}\nZ_{m_2}\right)= \Psi\left(\nZ_{m_1}\right)\Psi\left(\nZ_{m_2}\right) \\
&=\left\{-m_1, -(m_1-2), \ldots, m_1-2, m_1\right\}\left\{-m_2, -(m_2-2), \ldots, m_2-2, m_2\right\}.
\end{align*}
If we assume $m_1\leq m_2$, then by Theorem~\ref{thm:ford} we have
\begin{align}\label{eq:psi_m12}
\left|\Psi_n\left(\nZ_{m_1}\nZ_{m_2}\right)\right|\asymp\frac{n}{(\log m_1)^{\delta}(\log\log m_1)^{3/2}} 
\text{ as } \quad m_1 \rightarrow\infty. 
\end{align}
From Theorem~\ref{thm:basic_conditions} and (\ref{eq:psi_m12}) we get
\begin{align*}
\left|\Psi_n\left(\text{Comp}\left(\N[(\Z/2\Z)]\right)_n\right)\right|&=\left|\Psi_n\left(\bigcup_{\substack{d|n, \\ 1<d<n}}\nZ_d\nZ_{n/d}\right)\right|\\
&\leq 2\sum_{\substack{d | n\\ 1<d\leq \sqrt n }}\frac{n}{(\log d)^{\delta}(\log\log d)^{3/2}}.\\
&\ll_k \frac{n}{(\log P^-(n))^{\delta}(\log\log P^-(n))^{3/2}}.
\end{align*}
Considering only $\nZ_{P^-(n)}\nZ_{n/P^-(n)}$, we have 
\begin{align*}
\left|\Psi_n\left(\text{Comp}\left(\N[(\Z/2\Z)]\right)_n\right)\right|&\geq\left|\Psi_n\left(\nZ_{P^-(n)}\nZ_{n/P^-(n)}\right)\right|\\
&\gg \frac{n}{(\log P^-(n))^{\delta}(\log\log P^-(n))^{3/2}}.
\end{align*}
This proves the required result for $l=1$.

For $l>1$, instead of $\Psi$ we consider all the nontrivial characters of $\left(\Z/2\Z\right)^l$. Denote these characters by
\[\Psi^{(1)},\ldots,\Psi^{(2^l-1)}.\]
We can show that 
\[\Psi^{(j)}_n\left(\N[(\Z/2\Z)^l]_{n}\right)=\{-n, -n+2,\ldots,n-2, n\}.\]
Similar to the case $l=1$, using Theorem~\ref{thm:ford}, we can show that for $l> 1$ and $m_1\leq m_2$,
\[\Psi^{(j)}_n\left(\N[(\Z/2\Z)^l]_{m_1}\N[(\Z/2\Z)^l]_{m_2}\right)\asymp \frac{n}{(\log m_1)^{\delta}(\log\log m_1)^{3/2}}.\]
If we know $\Psi^{(j)}_n(a)$ for $j=1,\ldots,2^l-1$ and that $a\in \N[(\Z/2\Z)^l]_{n}$, then we can compute $a$. This gives
\begin{align*}
&\left|\text{Comp}\left(\N[(\Z/2\Z)^l]_n\right)\right|\leq \prod_{j=1}^{2^l-1}\left|\Psi^{(j)}_n\left(\text{Comp}\left(\N[(\Z/2\Z)^l]_{n}\right)\right)\right|\\
& \leq \prod_{j=1}^{2^l-1}\left(\sum_{\substack{d | n\\ 1<d\leq \sqrt n }}\frac{n}{(\log d)^{\delta}(\log\log d)^{3/2}}\right)
\ll \prod_{j=1}^{2^l-1} 2^k \frac{n}{(\log P^-(n))^{\delta}(\log\log P^-(n))^{3/2}}\\
&\ll_{k, l}\left(\frac{n}{(\log P^-(n))^{\delta}(\log\log P^-(n))^{3/2}}\right)^{2^l-1}.
\end{align*}
This completes the proof for $l> 1$.
\section*{Acknowledgement}
The author is supported by Grant 227768 of the Research Council of Norway, and this work was carried out when he was a research fellow at the Institute of Mathematical Sciences, Chennai.
The author would also like to thank Anirban Mukhopadhyay and Amritanshu Prasad for several fruitful discussions.


\begin{thebibliography}{3}
\bibitem{ford1}
Kevin Ford.
\newblock The distribution of integers with a divisor in a given interval.
\newblock {\em Ann. of Math.}, 168(2):367--433, 2008.

\bibitem{ford2}
Kevin Ford.
\newblock Integers with a divisor in {$(y,2y]$}.
\newblock {\em Anatomy of integers, CRM Proc. Lecture Notes}, 46:65--80, 2008.

\bibitem{golan}
Jonathan S. Golan.
\newblock {\em Semirings and their applications}.
\newblock Kluwer Academic Publishers, Dordrecht, 1999.

\bibitem{PHS}
G. Picci, J. M. van den Hof, J. H. van Schuppen
\newblock Primes in several classes of the positive matrices.
\newblock {\em Linear Algebra Appl.}, 277: 149--185, 1998.
%
%
%
\end{thebibliography}
\end{document}